\documentclass[11pt]{article}
\usepackage{amsmath,amssymb,a4wide}

\newtheorem{thm}{Theorem}
\newtheorem{co}[thm]{Corollary}
\newtheorem{lem}[thm]{Lemma}

\newtheorem{assumption}[thm]{Assumption}

\newtheorem{pr}[thm]{Proposition}

\newtheorem{assu1}[thm]{Assumption}

\newtheorem{definition}[thm]{Definition}

\newtheorem{example}[thm]{Example}

\newtheorem{remark}[thm]{Remark}

\newtheorem{protocol}[thm]{Protocol}

\newcommand{\dd}{\mathrm{d}}

\newcommand{\D}{\mathbb{D}}
\newcommand{\F}{\mathbb{F}}

\newcommand{\N}{\mathbb{N}}

\newcommand{\Z}{\mathbb{Z}}

\newcommand{\Section}[1]{\section{#1}}

\newcommand{\openbox}{\leavevmode
  \hbox to.77778em{%
    \hfil\vrule
  \vbox to.675em{\hrule width.6em\vfil\hrule}%
  \vrule\hfil}}
\newcommand{\proofname}{Proof}

\newenvironment{proof}[1][\proofname]{\par\normalfont
  \trivlist\item[\hskip\labelsep\itshape #1:]\ignorespaces
  }{\hspace*{1cm}\hspace*{\fill}\openbox \medskip\endtrivlist}

\newenvironment{proofofthm}[1][Proof of Proposition \ref{prop}]{\par\normalfont
  \trivlist\item[\hskip\labelsep\itshape #1:]\ignorespaces
  }{\hspace*{1cm}\hspace*{\fill}\openbox \medskip\endtrivlist}

\title{Natural Density of Rectangular Unimodular Integer Matrices
  \footnote{Partially supported by SNF grant No. 121874 and Armasuisse.}
}%
\date{\today}
\author{G\'erard Maze, Joachim Rosenthal and Urs Wagner\\
{\small {\em e-mail:\/} \{gmaze,rosen,urs.wagner\}@math.uzh.ch \vspace{-1mm} }\\
{\small Mathematics Institute\vspace{-1mm}}\\
{\small University of Z\"urich\vspace{-1mm}}\\
{\small Winterthurerstr 190, CH-8057 Z\"urich, Switzerland }
\vspace{3mm} }

\begin{document}\maketitle
\thispagestyle{empty}
\begin{abstract}
  In this paper, we compute the natural density of the set of $k
  \times n$ integer matrices that can be extended to an invertible $n
  \times n$ matrix over the integers. As a corollary, we find the
  density of rectangular matrices with Hermite normal form $\left[
    O_{k\times (n-k)} \, I_k \right]$. Connections with
  Ces\`aro's Theorem on the density of coprime integers and
  Quillen-Suslin's Theorem are also presented.
\end{abstract}

\vspace{3mm}
\vspace{3mm}

\Section{Introduction and Main Result}

Given a commutative ring $R$ with 1, the notion of invertible $n
\times n$ matrix is well defined, and can be characterized by the
condition that the determinant of such a matrix is a unit in $R$.
Given a $k \times n$ matrix $A$, the question of whether it can
be completed by an $(n-k) \times n$ matrix into an $n \times n$
invertible matrix over $R$ has raised several interesting
problems in the past. For instance, the celebrated Quillen-Suslin
Theorem \cite{quillen,suslin}, previously known under the name
Serre Conjecture, deals with the case when $R=\F[x_1,...,x_l]$,
with $\F$ a field and $k=1$. It can be shown that this case also
contains in essence the general case $1 \leq k \leq n$, see
\cite{youla}. The theorem states that over this ring, the
following three properties are equivalent
\cite{quillen,suslin,youla}:
\begin{enumerate}
\item  $A$ can be completed into an $n \times n$ invertible matrix,
\item there exists an  $n \times k$ matrix $B$ such that $AB= I_k$,
  where $I_{k}$ is the $k \times k$ identity matrix,
\item the $k \times k$ minors of $A$ have no common zeros.
\end{enumerate}

When the ring $R$ is a PID it is a direct consequence of the
Smith normal form (see e.g.~\cite{macduffee}) that Conditions 1.
and 2. are again equivalent and these conditions are equivalent
to the fact that the gcd of the $k\times k$ minors is equal to 1.

In the sequel we will adopt the usual convention and call a
$k\times n$ matrix $A$ over some ring $R$  {\it unimodular} as
soon as $A$ can be extended to $n\times n$ invertible matrix. 

Over the ring of integers various related results on unimodular
matrices are known in the literature. E.g. Zhan~\cite{zhan} showed
that any partial $n\times n$ matrix with $n$ given entries not lying
on the same row or column can be completed into a unimodular matrix.
Another result is due to Fang~\cite{fa} who showed that, if the
diagonal of a square matrix is let free, then it can be completed into
a unimodular one.

Our focus in this paper will be on the ``probability'' that a random
$k\times n$ integer matrix is unimodular. Related to our problem is a
classical result due to Ces\`aro~\cite{cesaro1,cesaro2,cesaro3} (see
also \cite{sylvester1,sylvester2} and the historical remarks below)
which states that the ``probability'' that two randomly chosen
integers are coprime is $\frac{1}{\zeta(2)}=\frac{6}{\pi^2}$, where
$\zeta$ denotes Riemann's zeta function. A re-statement of Ces\`aro's
result is then: the probability that a random $1\times 2$ integer
matrix is unimodular is $\zeta(2)^{-1}$.

In order to make the notion of probability precise we first remark
that the uniform distribution over the set $\Z^m$ has little
meaning. For this reason researchers often use the concept of {\em
  natural density} when stating probability results in $\Z^m$ or more
general infinite modules $R^m$. In the following we briefly explain
this concept. Let $S\subset \Z^m$ be a set. Define the upper
(respectively lower) natural density as
\[
\overline{\D}(S)=\limsup_{B \rightarrow \infty} \frac{|S \cap
\left[-B,B\right[^m|}{(2B)^m} \; \; , \;\;
\underline{\D}(S)=\liminf_{B \rightarrow \infty} \frac{|S \cap
\left[-B,B\right[^m|}{(2B)^m}.
\]
When both limits are equal one defines the natural density of the
set $S$ as:
\begin{equation}
  \D(S):=\overline{\D}(S)=\underline{\D}(S).
\end{equation}

The following properties of natural density are readily verified: If
$S^c$ denotes the complement of $S$ then $\D(S^c)=1-\D(S)$, whenever
$S$ has a well-defined density. Similarly if $\{S_i\}_{i\in I}$ is a
set of subsets of $\Z^m$ with well defined densities $\D(S_i)$ and if $S =
\cup_{i \in I} S_i$, then:
\begin{equation}
\overline{\D}(S) \leq \sum_{i \in I}\D(S_i).
\end{equation}
Readers interested in more background on the notion of natural
densities of sets of integer matrices are referred to~\cite{hetzel},
even though this paper attributes the result of Ces\`aro to Mertens,
something which was probably triggered by a remark in~\cite{hardy}.

In fact the exact fatherhood of the result of Ces\`aro appears to be
inexactly described in several occasions in the literature. The story
behind this historical misunderstanding seems to be the following. The
problem of evaluating the probability that two random integers are
coprime appears in a 1881 question raised by Ces\`aro
\cite{cesaro1}. Two years later, Sylvester \cite{sylvester1} and
Ces\`aro \cite{cesaro2} independently publish their
solutions. Interestingly, two proofs are presented in
\cite{sylvester1}: Sylvester's own argument is based on Farey series, and
a more ``probabilistic'' argument of Franklin is also presented with
his permission. In the footnote of a 1888 paper \cite{sylvester2},
Sylvester publishes a similar proof, and mentions that Ces\`aro
``claimed the prior publication`` of the result. Sylvester's argument
is based on Farey series, and the remark in~\cite{hardy} makes a
connection between an earlier work of Mertens~\cite{mertens} of 1874
on the average value of Euler $\varphi$ function and Farey series,
which probably triggered the remark in~\cite{hetzel} but at no point
in~\cite{hardy} the result is attributed to Mertens. If we want to
associate the aforementioned probability with the average value of
$\varphi$, then it is legitimate to go back to an 1849 paper of
Dirichlet \cite{dirichlet} where the value $\frac{6}{\pi^2}$ appears
for the first time. The case of $k$ coprime integers, $k>2$, is also
presented for first time by Ces\`aro in 1884 \cite{cesaro3}. The
result is rediscovered in 1900 by D.N.Lehmer \cite{lehmer}, apparently
independently (see also \cite{nymann}).

The major result of our paper is a matrix version of Ces\`aro's
theorem:

\begin{pr}  \label{prop}
  Let $1 \leq k < n$ be integers. The natural density $\dd_{k,n}$ of
  $k \times n$ unimodular matrices with integer entries is given by
  $$
  \dd_{k,n} = \left( \prod_{j=n-k+1}^{n} \zeta(j) \right)^{-1}.
  $$
\end{pr}
In other words, the ``probability'' that a ``random'' $k \times n$
integer matrix can be extended into a matrix in $\mathrm{GL}_n(\Z)$ is
given by $\left( \zeta(n)\cdot \ldots \cdot \zeta(n-k+1)
\right)^{-1}$. Since $\zeta(n)$ converges rapidly towards 1, if
$k=n-d$, the above density converges rather fast to the limit
$\dd_{d}$ with
\[
\dd_1 = \left( \prod_{j=2}^{\infty} \zeta(j) \right)^{-1}
=0.43575707677... \;\; \mbox{ and } \;\; \dd_2 = \zeta(2) \cdot \dd_1
, \; \dd_3 = \zeta(2) \cdot \zeta(3) \cdot \dd_1 \ldots
\]
When $n$ is not too small, say $n>4$, $\dd_1$ is a good approximation
of the proportion of integer matrices that can be completed by a row
into an invertible matrix. We would like to point out that the proof
of the above proposition is independent from the result of Ces\`aro
and as such a new proof of his theorem is given when the above
proposition is considered with $k=1$. Proposition~\ref{prop} above
gives a probabilistic extension of the simplest case of the
Quillen-Suslin Theorem. The concept of natural density does not extend
naturally to the ring $\F[x_1,\ldots,x_n]$ and thus the existence of a
direct extension of our result to the general case is unclear. In the
next section we will prove Proposition~\ref{prop}.

\section{Proof of the Main Result}

Let us fix some notations for the rest of the article. Let $1 \leq t
\leq k < n$ be integers. Given a $k \times n$ matrix $A$ with
integer entries, the determinants of the $\binom{n}{t} \binom{k}{t}$
$t \times t$ submatrices of $A$ formed by the intersection of any
subset of $t$ column and $t$ row vectors of $A$ are called the
$t$-minors of $A$. When $t=k$, they are called the full rank minors of
$A$. The set of primes is denoted by $\mathbb{P}$ and when no
confusion is possible, the set $\Z^{k\times n}$ will be identified
with $\Z^{kn}$.\\

The strategy of the proof of the above proposition is to localize the
computation of the density at every prime, and then lift the
information up in order to extract the exact density over $\Z^{kn}$.
The proof of the above proposition relies on the next lemmas. As noted
before, the PID version of Quillen-Suslin's Theorem gives directly
the following lemma:

\begin{lem}                       \label{lem1}
  A $k \times n$ matrix $A$ with integer entries is unimodular if and
  only if the full-rank minors of $A$ are coprime.
\end{lem}

Over a field $\F$, the rank of a $k \times n$ matrix, i.e., the
maximal number of rows that are linearly independent over $\F$, is
equal to the largest integer $t$ such that there exists a non-zero $t
\times t$ minor. A classical computation shows that over the prime
finite field $\F=\Z/p\Z$, if $F_p$ is the set of full rank $k \times
n$ matrices, its cardinality satisfies $|F_p| = \prod_{j=0}^{k-1}
(p^n-p^j)$. See \cite{lidl} for the details. This can be extended as
follows:

\begin{lem}                             \label{lem2}
  Let $S$ be a finite set of prime numbers. The density of $k \times n$
  matrices with integer entries for which the $\gcd$ of the full rank
  minors are coprime to all primes in $S$ is given by
  \[
  \prod_{j=n-k+1}^{n} \prod_{p \in S} \left( 1-\frac{1}{p^j} \right).
  \]
\end{lem}

\begin{proof}
  Let us call $E_S \subset \Z^{kn}$ the set of $k \times n$
  matrices for which the gcd of the full rank minors are coprime
  to all primes in $S$. Let $N=\prod_{p \in S} p$ and recall the
  Chinese remainder theorem $\left( \Z / N \Z \right)\cong
  \prod_{p \in S} \left( \Z / p \Z \right)$.  Let $B$ be an
  integer that will go to infinity in the sequel. Write $B=qN+r$,
  with $q,r \in \N$, $0 \leq r < N$ and consider the map $\phi$
  obtained as the composition of maps
  \[
  \lbrack -qN,qN \lbrack^{kn} \longrightarrow \left( \Z / N \Z \right)^{kn}
  \longrightarrow  \prod_{p \in S} \left( \Z / p \Z \right)^{kn}
  \]
  where the first map is the quotient modulo $N$ and the second is the
  induced homomorphism given by the Chinese remainder theorem. Because
  of the above remark we have
  \[
  \phi(E_S \cap \lbrack  -qN,qN \lbrack^{kn} ) = \prod_{p \in S} F_p.
  \]
  The first map is a $(2q)^{kn}$-to-$1$ map, i.e., each fiber contains
  exactly $(2q)^{kn}$ elements, and the second map is an isomorphism. Thus
  \[
  |E_S \cap \lbrack-qN,qN\lbrack^{kn}| = (2q)^{kn} \cdot  \prod_{p \in S} |F_p|.
  \]
  We have the disjoint union $\lbrack-B,B\lbrack^{kn} =
  \lbrack-qN,qN\lbrack^{kn} \sqcup \left( \lbrack-B,B\lbrack^{kn}
  \setminus \lbrack-qN,qN\lbrack^{kn} \right)$. The difference of
  hypercubes has a volume bounded by the volume of $2kn$
  hyper-rectangles of side area $B^{kn-1}$ and height $r$ which gives
  $0 \leq | \lbrack-B,B\lbrack^{kn} \setminus
  \lbrack-qN,qN\lbrack^{kn} | < 2knrB^{kn-1}$ and therefore we have
  \[
  |E_S \cap \lbrack-B,B \lbrack^{kn}| = |E_S \cap \lbrack-qN,qN \lbrack^{kn}| + \rho
  \]
  with $0 \leq \rho < 2knrB^{kn-1}$. Thus
  \begin{eqnarray*}
    \frac{|E_S \cap [-B,B[^{kn}|}{(2B)^{kn}} & = & \frac{ |E_S \cap
          [-qN,qN[^{kn}| }{(2qN)^{kn}} \frac{(2qN)^{kn}}{(2B)^{kn}} +
            \frac{\rho}{(2B)^{kn}}\\ 
            & = &  \prod_{p \in S} \frac{|F_p|}{p^{nk}} \cdot
            \left(1-\frac{r}{B} \right)^{kn}+ O(1/B).
  \end{eqnarray*}
  Finally, we have 
  \[
  \D(E_S) = \lim_{B \rightarrow \infty}  \prod_{p \in S} \frac{|F_p|}{p^{nk}} \cdot
  \left(1-\frac{r}{B} \right)^{kn}+ O(1/B) = 
  \prod_{p \in S} \frac{|F_p|}{p^{nk}} = \prod_{j=n-k+1}^{n} \prod_{p
    \in S} \left( 1-\frac{1}{p^j} \right).
  \]
\end{proof}

Extending the definition of $E_S$, let us call $E_t \subset \Z^{kn}$
the set of $k \times n$ matrices for which the gcd of the full rank
minors are coprime with the first $t$ primes $2,3,\ldots,p_t$. Note
that the sequence $E_t$ is a decreasing sequence of sets, i.e., $E_i
\subset E_j$ if $i \geq j$ and if $E = \cap_{t\in \N} E_t$, then Lemma
\ref{lem1} implies that $E$ is the set of $k \times n$ unimodular
matrices with integer entries. In order to prove Proposition
\ref{prop}, we have to prove that $\D(E)$ exists and compute its value
$\dd_{k,n} = \D(E)$. Since we know $\D(E_t)$ for all $t$, it is
tempting to prove the proposition by simply letting $t$ going to
infinity in the expression given by Lemma \ref{lem2}, but this is an
invalid argument in general. Indeed the example of $E_t = \lbrack t,
\infty \lbrack$ with $E = \emptyset$ shows that it is possible to have
a sequence of sets $E_i \subset E_j$ if $i \geq j$ with $E =
\cap_{t\in \N} E_t$ and $\D(E)=0 \neq 1 = \lim_{t \rightarrow \infty}
\D(E_t)$ since $\D(E_t) = 1$, $\forall t$. The next lemma describes
how to avoid this pathological case. Let us recall that for any real
sequences $a_n,b_n$, $\liminf a_n + \liminf b_n \leq \liminf (a_n +
b_n)$ and $\limsup (a_n + b_n) \leq \limsup a_n + \limsup b_n$ and
$\limsup -a_n = - \liminf a_n$.

\begin{lem}\label{lem3}
  Let $E_t$ be a sequence of decreasing sets in $\Z^m$ such that
  $\D(E_t)$ exists for all $t$ and converges to $\dd$. Let $E =
  \cap_{t\in \N} E_t$. If $\lim_{t \rightarrow \infty}
  \overline{\D}(E_t \setminus E) = 0$, then $\D(E)$ exists and is
  equal to $\dd$.
\end{lem}

\begin{proof}
  We will use the disjoint union $E_t = E \sqcup (E_t \setminus
  E)$. Since
  \[
  |E_t \cap \lbrack-B,B \lbrack^m| = |E \cap
    \lbrack-B,B \lbrack^m| +|(E_t \setminus E) \cap
    \lbrack-B,B \lbrack^m|,
  \]
  we have
  \[
  \frac{|E \cap \lbrack-B,B \lbrack^m|}{(2B)^m} = \frac{|E_t \cap
    \lbrack-B,B \lbrack^m|}{(2B)^m} + \frac{-|(E_t \setminus E) \cap
    \lbrack-B,B \lbrack^m|}{(2B)^m}.
  \]
  Taking the $\liminf$ and the $\limsup$, and since $\lim_{B
    \rightarrow \infty} \frac{|E_t \cap \lbrack-B,B
    \lbrack^m|}{(2B)^m}=\D(E_t)$, we have
  \[
  \D(E_t) - \overline{\D}(E_t \setminus E) =  \underline{\D}(E) \leq 
  \overline{\D}(E) = \D(E_t) - \underline{\D}(E_t \setminus E).
  \]
  The result follows when $t \rightarrow \infty$ since $0 \leq 
  \underline{\D}(E_t \setminus E) \leq \overline{\D}(E_t \setminus E)
  \longrightarrow 0$.
\end{proof}

\begin{proofofthm}
  We will use the previous lemma. The set $E_t \setminus E$ is the set
  of $k \times n$ matrices $A$ for which there exists a prime $p$ with
  $p>p_t$ so that $p$ divides the gcd of the full rank minors of
  $A$. For each prime $p$, let us define $H_p$ to be the set of $k
  \times n$ matrices whose gcd of the full rank minors is divisible by
  $p$. Then $ E_t \setminus E = \cup_{p > p_t} H_p$. Note that Lemma
  \ref{lem2} applied to $S=\{p\}$ implies that $H_p$ has a density
  equal to
  \[
  \D(H_p) = 1-\D(H_p^c) = 1-\prod_{j=n-k+1}^{n} \left(1-\frac{1}{p^j}\right).
  \]
  By induction on the number of factors, one readily verifies that for
  real numbers $0<x_j<1$, we have the inequality $\prod_{j=0}^{n}
  (1-x_j) > 1- \sum_{j=0}^{n} x_j$, which applied to the above product
  gives
  \[
  \D(H_p) = 1-\prod_{j=n-k+1}^{n} \left(1-\frac{1}{p^j}\right) < \sum_{j=n-k+1}^{n}
  \frac{1}{p^j} < \frac{1}{p^{n-k}(p-1)} < \frac{2}{p^2}.
  \]
  Finally, we have
  \[
  \overline{\D}(E_t \setminus E) = \overline{\D}(\cup_{p > p_t} H_p)
  \leq \sum_{p > p_t} \D(H_p) < \sum_{p > p_t} \frac{2}{p^2}
  \]
  which shows that $\lim_{t \rightarrow \infty} \overline{\D}(E_t
  \setminus E) = 0$ since the last series is the tail of the
  convergent series $ \sum_{p \in \mathbb{P}} \frac{2}{p^2}$. We can
  therefore apply the previous lemma and conclude that
  $\dd_{k,n}=\D(E)$ exists and is equal to
  \[
  \dd_{k,n}=\D(E) = \prod_{j=n-k+1}^{n} \prod_{p \in \mathbb{P}}
  \left(1-\frac{1}{p^j}\right) = \prod_{j=n-k+1}^{n} \zeta(j)^{-1}.
  \]

\end{proofofthm}

\section{Concluding Remarks and Further Results}

Proposition~\ref{prop}  does not cover the square
case $k=n$. Indeed, the zeta function is not defined when $n-k+1=1$
since $\zeta$ has a pole of order 1 at $x=1$. The following
result covers the square case.
\begin{lem}
The natural density of $n\times n$ unimodular matrices is
$\dd_{n,n}=0$. 
\end{lem}
Before we give a proof we remark that the product formula in
Proposition~\ref{prop} naturally has an extension to zero as 
 $\lim_{x\rightarrow 1} (\zeta(x))^{-1}=0$. 
\begin{proof}
  Each $n \times n$ matrix with $n^2-1$ entries in the range
  $\lbrack-B,B \lbrack$ can be completed by at most two values in
  order for this matrix to be unimodular, due to the Lagrange
  expansion of the determinant, which must be $\pm 1$. As such
  there are at most $2(2B)^{n^2-1}$ unimodular matrices with
  entries in $\lbrack-B,B \lbrack$.  The conclusion follows since
  $\dd_{n,n} \leq \lim_{B \rightarrow \infty} 2(2B)^{n(n-1)} /
  (2B)^{n^2} = 0$.
\end{proof}

The result of Proposition~\ref{prop} can also be used in the
determination of the natural density of $k \times n$ matrices whose
Hermite normal form (HNF) is very simple. Recall that the HNF of a $k
\times n$ matrix $A$ is the unique $k \times n$ matrix $H$ of the
following form 
\[
\left(
\begin{array}{cccccccc}
0 & 0 & \ldots & 0 & h_1 & h_{1,2} & \ldots & h_{1,n-k} \\
0 & 0 & \ldots & 0 & 0 & h_2 &  \ldots & h_{2,n-k}\\
\vdots & \vdots & \ddots & \vdots & \vdots &  \ddots &  \ddots &
\vdots\\
0 & 0 & \ldots & 0 & 0 & \ldots & 0 & h_{k}
\end{array}
\right)
\]
where $0 \leq h_{i,j} < h_i$, such that there exists $U \in
\mbox{Gl}_n (\Z)$ with $A=HU$, see e.g. \cite{cohen}. We have then the
following result.

\begin{thm}
  The density of $k \times n$ matrices whose Hermite normal form is
  the block matrix $\lbrack O_{k\times (n-k)} \, I_k \rbrack$ is
  $\dd_{k,n}$.
\end{thm}
 
\begin{proof}
The result follows from the fact that the set $E$ of $k \times n$
unimodular matrices over $\Z$ is the equal to the set $L$ of $k \times
n$ matrices whose HNF is the block matrix $\lbrack O_{k\times (n-k)}
\, I_k \rbrack$. Let us show that $E \subset L$. If $A \in E$
then there exists a unimodular square matrix $B$ such that $A$
consists in the last $k$ rows of $B$. The HNF of $B$ is the identity
matrix $I_n$. Indeed, it is a upper diagonal square matrix with
integer entries whose determinant is 1, which forces the diagonal
entries to be 1. The size condition of the coefficient of the HNF
forces the entries above the diagonal to be 0. Thus we have
$BU^{-1}=I_n$ for some $U \in \mbox{Gl}_n (\Z)$, which gives
$AU^{-1}=\lbrack O_{k\times (n-k)} \, I_k \rbrack$. The
uniqueness of the HNF shows that $A \in L$. Let us show now that $L
\subset E$ by showing that the gcd of the $k \times k$ minors of a
matrix in $L$ are coprime. The equation $A = \lbrack O_{k\times (n-k)}
\, I_k \rbrack U$ shows that the $k \times k$ minors of $A$
are equal to the $k \times k$ minors of $\lbrack O_{k\times (n-k)} \,
I_k \rbrack$ which are 0 or 1, and thus coprime. This finishes
the proof of the corollary.
\end{proof}

Taking into account that the Smith normal form of a matrix
\cite{cohen} is obtained from the Hermite normal form via row
operations, an immediate consequence of this theorem is:

\begin{co}
The density of $k \times n$ matrices whose Smith normal form is the
block matrix $\lbrack O_{k\times (n-k)} \, I_k \rbrack$ is
$\dd_{k,n}$.
\end{co}


\begin{thebibliography}{10}

\bibitem{cesaro1}
         Ces\`aro, E.
         \newblock Question propos\'ee 75.
         \newblock {\it Mathesis}, 1 (1881), p. 184.

\bibitem{cesaro2}
         Ces\`aro, E.
         \newblock Question 75 (Solution).
         \newblock {\it Mathesis}, 3 (1883), pp. 224--225.


\bibitem{cesaro3}
        Ces\`aro, E.
        \newblock Probabilit\'e de certains faits arithm\'ethiques.
        \newblock {\it Mathesis} 4 (1884), pp. 150--151.

\bibitem{cohen}
        Cohen, H.
        \newblock {\it A course in computational algebraic number theory},
        Springer Graduate Texts In Mathematics, 1995.

\bibitem{dirichlet}
        Dirichlet, P. G. L.
        \newblock \"Uber die bestimmung der mittleren Werthe in
          der Zahlentheorie
        \newblock {\it Abhand. Ak. Wiss. Berlin} (1849),
        pp. 63--83. Reprinted in {\it Werke}, Vol. 2, pp. 49--66. 


\bibitem{fa}
        Fang, M.
        \newblock On the completion of a partial integral matrix to a 
        unimodular matrix.
        \newblock {\it Linear Algebra and its Applications}
        Volume 422, Issue 1, 1 April 2007, pp. 291--294.


\bibitem{hardy}
        Hardy, G. H., Wright, E. M.
        \newblock {\it An Introduction to the Theory of Numbers},
        \newblock 5d ed. Oxford University Press, 1979.

\bibitem{hetzel}
        Hetzel, A. J., Liew, J. S. and Morrison, K. E.
        \newblock The probability that a matrix of integers is diagonalizable.
        \newblock {\it American Mathematical Monthly}, Volume 114,
        Number 6, June-July 2007, pp. 491--499. 

\bibitem{lehmer}
        Lehmer, D.N.
        \newblock  Asymptotic evaluation of certain totient sums.
        \newblock {\it American Journal of Mathematics}, Vol. 22,
        No. 4 (Oct., 1900), pp. 293-335. 

\bibitem{lidl}
        Lidl, R., Niederreiter, H.
        \newblock {\em Finite fields},
        \newblock  Second edition, Cambridge. University Press, 1997.

\bibitem{macduffee}
        MacDuffie, C. C.
        \newblock {\em The Theory of Matrices},
        \newblock Chelsea Publ. Co., New York, 1946.

\bibitem{mertens}
        Mertens, F.
        \newblock \"Uber einige asymptotische Gesetze der
        Zahlentheorie.
        \newblock {\it " J. reine angew. Math.} 77, pp. 46--62, 1874.

\bibitem{nymann}
        Nymann, J. E. 
        \newblock On the probability that $k$ positive integers are
        relatively prime. 
        \newblock {\it J. Number Th.}, 7 (1975), pp. 406-412.

\bibitem{quillen}
        Quillen, D.
        \newblock Projective modules over polynomial rings.
        \newblock {\it Inventiones Mathematicae}, No. 36, 1976,
        pp. 167--171

\bibitem{suslin}
        Suslin, A. A.
        \newblock Projective modules over polynomial rings are free.
        \newblock {\it Soviets Mathematics}, No. 17, Volume 4, 1976,
        pp. 1160--1164. 

\bibitem{sylvester1}
        Sylvester, J.J.
        \newblock Sur le nombre de fractions ordinaires in\'egales
        qu'on peut exprimer en se servant de chiffres qui n'exc\`ede
        pas un nombre donn\'e.
        \newblock {\it Comptes Rendus de l'Acad. Sci. Paris}, XCVI
        (1883), pp. 409--413.
        Reprinted in  Baker, H. F. ed.,
        \newblock {\it The collected mathematical papers of James
          Joseph Sylvester}, Volume 4, Cambridge University press, p. 86.

\bibitem{sylvester2}
        Sylvester, J.J.
        \newblock  On certain inequalities relating to prime numbers.
        \newblock {\it Nature}, No. 38, 259-262, 12 July 1888.

\bibitem{youla}
        Youla D. C., Pickel P. F. 
        \newblock The Quillen-Suslin theorem and the structure of 
        $n$-dimensional elementary polynomial matrices.
        \newblock {\it IEEE transactions on circuits and systems},
        1984, Volume 31, Number 6, pp. 513--518.

\bibitem{zhan}
        Zhan, X.
        \newblock Completion of a partial integral matrix to a 
        unimodular matrix.
        \newblock {\it Linear Algebra and its Applications}, Volume
        414, Issue 1, 1 April 2006, pp. 373--377.



\end{thebibliography}
\end{document}